# UNITING CONTROL FOR FINITE TIME STABILIZATION OF NONLINEAR DETECTABLE SYSTEMS

Denis Efimov, Alexander L. Fradkov

*Abstract*—The stabilization of nonlinear systems under zero-state-detectability assumption or its analogues is considered. The proposed supervisory control provides a finite time practical stabilization of output and it is based on uniting local and global controllers. The global control ensures boundedness of solutions and output convergence to zero, while local one ensures finite time convergence to a predefined set into the zero dynamics set. Computer simulation illustrates advantages of the proposed algorithm.

I. INTRODUCTION

In many applications it is possible to design a global stabilizing control law (ensuring global boundedness of system solutions in the presence of disturbances) and local control law (guaranteeing optimality of the solutions in some sense without disturbances). In such cases it is desirable to design a united controller, which inherits properties of both local and global ones. For applications the quality of transient processes is very important and a global solution should be proposed with optimal properties of the local control law. Uniting control can be considered as a such solution if it coincides with the optimal controller near the desired set and provides boundedness of the system solution under disturbances.

The first solution to the uniting control problem was suggested in [32] where a dynamic time invariant control law was proposed converging to locally optimal control near the origin under special conditions only. In the paper [19] a static time invariant control was presented under condition of existence of a continuous path between global and local controls which is hard to verify. In the works [21], [23] an example was found, that does not allow any continuous or even discontinuous time invariant controls. Additionally, in these works several solutions of uniting control problem were proposed (continuous, discontinuous, hybrid and time-varying) for the case of vanishing external disturbances. A kind of uniting control for chained nonholonomic systems was developed in [22], where robust properties of such control law with respect to sufficiently small disturbances were analyzed. Uniting control under acting disturbances was considered in [6], [14]. Despite this success an evaluation of quality improvement achieved in closed by uniting control system has not been presented yet.

In this paper we are going to concentrate our attention on evaluation and comparison of *transient time* for systems governed by global and uniting controls. This problem becomes very important in some situations requiring zero-state-detectability assumption [4], [17] or its analogues [3], [24]−[27] for a desired set stabilization. Roughly speaking, in these cases it is possible to design global control laws ensuring robust stability and convergence of all trajectories to an invariant set $\Gamma$, into the set $\Gamma$ all trajectories reach for a desired subset $\gamma \subset \Gamma$ in some time. It will be shown that in some cases such systems can possess well defined time of convergence to the invariant set $\Gamma$ without guaranteed time of reaching for the desired subset $\gamma$. To compensate this shortage an uniting control is proposed. While global control provides global convergence and robustness, the local one ensures required time estimates for the desired set $\gamma$ reaching.

In the next section a motivating example is considered explaining the importance of transient time issue in the systems with detectability property. In section 3 an uniting control is designed solving the problem. In section 4 advances of the proposed solution are demonstrated via analytical design and computer simulation.

## II. MOTIVATING EXAMPLE

Let us consider a well known problem [1], [3], [12], [18] of swinging up a nonlinear pendulum:

$$\dot{x}_1 = x_2;$$
$$\dot{x}_2 = -\omega^2 \sin(x_1) + \cos(x_1)u,$$
(1)

where $x_1 \in R$ is angular position; $x_2 \in R$ is angular velocity; $\mathbf{x} = [x_1 \ x_2] \in R^2$ is state vector of the system; $\omega \in R$ is natural frequency of (1); $u \in [-u_m \ u_m]$ is control input, $0 < u_m < +\infty$. It is required to stabilize the upper unstable equilibrium of the pendulum (1) with coordinates $(n\pi, 0)$ for some $n = \pm 1, \pm 3, \ldots$. To this purpose let us utilize the energy based approach [12], [13], [24]:

$$u = -u_m \, sign\big([H(\mathbf{x}) - H^*]x_2 \cos(x_1)\big),$$
(2)

$$sign(u) = \begin{cases} 1 & if \ u \geq 0; \\ -1 & if \ u < 0, \end{cases}$$

where $H(\mathbf{x}) = 0.5 x_2^2 + \omega^2(1 - \cos(x_1))$ is the energy function of the pendulum (1), $H^* = 2\omega^2$ is the stabilized value of the energy $H$ corresponding to the upper equilibrium. As it was proven in [28] for $H = H^*$ system (1), (2) is detectable with respect to the upper equilibrium (below in the text the statement "A $\Rightarrow$ B" means that "if A is true, then B is satisfied"):

$$H(\mathbf{x}(t)) \equiv H^*, \ t \geq 0 \ \Rightarrow \ \lim_{t \to +\infty} \mathbf{x}(t) = (n\pi, 0), \ n = \pm 1, \pm 3, \ldots.$$
(3)

Unfortunately relation (3) does not provide any estimate of transient processes time length in the system. Denote as

$$T_H(\mathbf{x}_0) = \arg\inf_{t \geq 0} \{ |H(\mathbf{x}(\tau, \mathbf{x}_0)) - H^*| \leq 0.05 H^*, \forall \tau \geq t \},$$

$$T_x(\mathbf{x}_0) = \arg\inf_{t \geq 0} \{ |x_1(\tau, \mathbf{x}_0) - n\pi| \leq 0.05 \pi, \forall \tau \geq t \}$$

the time of reaching for 5% zone of the desired values for variables $H$ and $x_1$ correspondingly. Here $\mathbf{x}_0 \in R^2$ is initial conditions vector. Any of functions $T_H$ or $T_x$ can be chosen for estimation of transient time in the system. Time $T_H$ can be assigned and evaluated from control (2) since for the system (1), (2) the relation

$$\frac{d[H(\mathbf{x}(t)) - H^*]}{dt} = -u_m \ |x_2(t)\cos(x_1(t))| \times$$
$$\times sign\big(H(\mathbf{x}(t)) - H^*\big),$$

holds and for persistently excited signal $x_2(t)\cos(x_1(t))$ (see Appendix) the system possesses an upper estimate of convergence time. Unfortunately evaluation of time $T_x$ is more sophisticated since the property (3) does not impose any restrictions on the rate of $\mathbf{x}$ convergence. To illustrate peculiarities of time $T_x$ let us simulate the solutions of the system (1), (2) for initial conditions $-10 \leq x_1^0 \leq 10$, $-3 \leq x_2^0 \leq 3$ with values of parameters $\omega = 1$, $u_m = 0.1$. The results of simulation are shown in Fig. 1 (for the time interval $0 \leq t \leq 200$ sec). As it is possible to conclude from Fig. 1, the time $T_H$ is proportional to the initial deviation $H(\mathbf{x}_0) - H^*$ of the trajectories, while the time $T_x$ has chaotic nature due to (3).

According to Fig 1, there exist closely located initial conditions with seriously deviated values of time $T_x$.

There exists another variant of transient time evaluation based on probability approach suitable in some areas of application [11]. Denote

$$T_x(t) = \int_{\mathbf{x}_0 \in \Omega} \varphi_x(\mathbf{x}(t,\mathbf{x}_0))\,d\mathbf{x}_0\,,\quad \varphi_x(\mathbf{x}) = \begin{cases} 1 & \text{if } |\mathbf{x} - n\pi| \le 0.05\,\pi; \\ 0 & \text{otherwise}; \end{cases}$$

$$T_H(t) = \int_{\mathbf{x}_0 \in \Omega} \varphi_H(\mathbf{x}(t,\mathbf{x}_0))\,d\mathbf{x}_0\,,\quad \varphi_H(\mathbf{x}) = \begin{cases} 1 & \text{if } |H(\mathbf{x}) - H^*| \le 0.05\,H^*; \\ 0 & \text{otherwise}, \end{cases}$$

as the cumulative distributions of the variables $\mathbf{x}$ and $H$ with respect to 5% zone of the desired values. In this case

$$p_x(t) = \frac{d\,T_x(t)}{d\,t},\quad \int_0^{+\infty} p_x(t)\,dt = 1;\quad p_H(t) = \frac{d\,T_H(t)}{d\,t},\quad \int_0^{+\infty} p_H(t)\,dt = 1$$

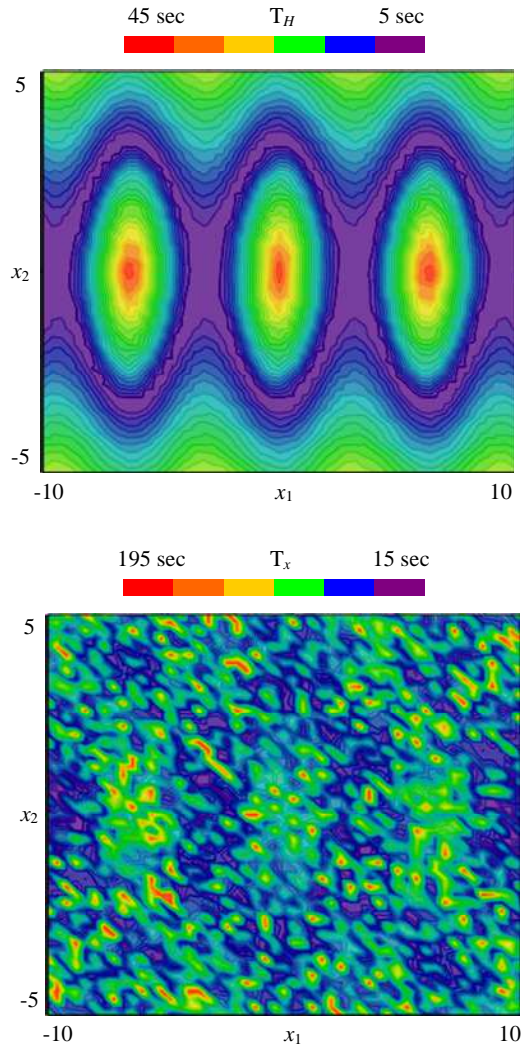

Fig. 1. Times $T_x$ and $T_H$ for system (1), (2).

are probabilities of appearance of a trajectory at time instant $t \ge 0$ in the corresponding 5% zone, and

$$E_x = -\int_0^{+\infty} p_x(t)\ln(p_x(t))\,dt\,,\quad E_H = -\int_0^{+\infty} p_H(t)\ln(p_H(t))\,dt$$

are the entropies of the system (1), (2) with respect to the 5% zones. The functions $T_x$, $T_H$, $p_x$, $p_H$ and the values $E_x$, $E_H$ can be considered as indirect characteristics of the transient time in the system (minimization of the entropies is equivalent to minimization of $T_H$ or $T_x$ dispersion).

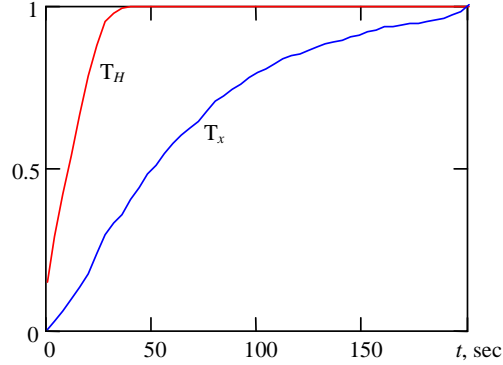

Fig. 2. Plots of distributions $T_x$ and $T_H$ for system (1), (2).

For the mentioned set of initial conditions and values of the system parameters (for the time interval $0 \leq t \leq 200$ sec) the plots of the functions $T_x$, $T_H$ are presented in Fig. 2 and

$$E_x = 9.2, \quad E_H = 6.6.$$

To calculate these values the interval $0 \leq t \leq 200$ was divided into $2 \times 10^4$ points. As it is possible to conclude from Fig. 2, after time $t_H = 50$ sec the system reaches the 5% zone for the energy $H$ with probability 0.99, while for the variable $\mathbf{x}$ such time is almost four times bigger ($t_x = 195$ sec). The same for the values of entropy (to compare the values of $E_x$, $E_H$ note, that entropy of the uniformly distributed process on $2 \times 10^4$ points is equal to $\ln(2 \times 10^4) = 9.9$).

Therefore, for this system under simple detectability assumption (3) there is a serious problem with the time of the upper equilibrium reaching. This problem becomes important in situations when it is necessary to provide a finite time convergence to a neighborhood of predefined position.

To point out a possible way of the problem solution let us note that for the linear model of pendulum

$$\begin{aligned} \dot{x}_1 &= x_2; \\ \dot{x}_2 &= -\omega^2 x_1 + u, \end{aligned} \quad (4)$$

with linear control

$$u = -K_1(x_1 - n\pi) - K_2 x_2 - \omega^2 x_1 \quad (5)$$

there is no such problem with transient time. And of course the problem described above is originated by nonlinearity of the system (1), (2). Though nonlinearities can not be neglected globally, the nonlinear system (1), (2) can be reduced to a linear one like (4), (5) locally, i.e. near the upper equilibrium. Combining the estimates for time needed to reach the neighborhood of the upper equilibrium and local time estimates for linearized system near the equilibrium it is possible to provide desired finite time convergence to a predefined subset near the equilibrium.

The problem can be expressed as the problem of stabilization with respect to two outputs, the set of zeros for the first one

$y = H(\mathbf{x}) - H^* = 0.5 x_2^2 + \omega^2(1 - \cos(x_1)) - 2\omega^2$ indicates the desired energy level, and the set of zeros of the output $\psi = [1 + \cos(x_1), x_2]$ corresponds to stabilized equilibriums. Thus, the idea of the paper consists in uniting local and global controllers designed for different outputs stabilization under detectability assumption (which establishes a relation among the outputs) with guaranteed length of transient processes.

## III. Main Result

Consider the nonlinear system

$$\dot{\mathbf{x}} = \mathbf{f}(\mathbf{x}, \mathbf{u}, \mathbf{d}), \quad \mathbf{y} = \mathbf{h}(\mathbf{x}), \quad \psi = \eta(\mathbf{x}), \tag{6}$$

where $\mathbf{x} \in R^n$ is the state vector; $\mathbf{u} \in R^m$ is the control input vector; $\mathbf{d} \in R^l$ is the vector of disturbances; $\mathbf{y} \in R^p$, $\psi \in R^q$ are auxiliary outputs; functions $\mathbf{f}: R^{n+m+l} \to R^n$, $\mathbf{h}: R^n \to R^p$ and $\eta: R^n \to R^q$ are continuous and locally Lipschitz. Euclidean norm will be denoted as $|\mathbf{x}|$, and $\|\mathbf{u}\|_{[t_0,t]}$ denotes the $L_\infty^m$ norm of the input (inputs $\mathbf{u}(t)$, $\mathbf{d}(t)$ are measurable and locally essentially bounded functions $\mathbf{u}: R_+ \to R^m$, $\mathbf{d}: R_+ \to R^l$, $R_+ = \{\tau \in R : \tau \geq 0\}$):

$$\|\mathbf{u}\|_{[t_0,T]} = \mathrm{ess\,sup}\{|\mathbf{u}(t)|, t \in [t_0, T]\},$$

if $T = +\infty$ then we will simply write $\|\mathbf{u}\|$. We will denote as $\mathcal{M}_{R^m}$ the set of all such Lebesgue measurable inputs $\mathbf{u}$ with property $\|\mathbf{u}\| < +\infty$. Denote as $\mathbf{x}(t, \mathbf{x}_0, \mathbf{u}, \mathbf{d})$ the solution of the system (6) with initial conditions $\mathbf{x}_0 \in R^n$ for inputs $\mathbf{u} \in \mathcal{M}_{R^m}$ and $\mathbf{d} \in \mathcal{M}_{R^l}$, $t \geq 0$, $\mathbf{y}(t, \mathbf{x}_0, \mathbf{u}, \mathbf{d}) = \mathbf{h}(\mathbf{x}(t, \mathbf{x}_0, \mathbf{u}, \mathbf{d}))$, $\psi(t, \mathbf{x}_0, \mathbf{u}, \mathbf{d}) = \eta(\mathbf{x}(t, \mathbf{x}_0, \mathbf{u}, \mathbf{d}))$ (we will simple write $\mathbf{x}(t)$, $\mathbf{y}(t)$ or $\psi(t)$ if all other arguments are clear from the context). The solutions are defined on some finite interval $[0, T)$; if $T = +\infty$ for every initial state $\mathbf{x}_0 \in R^n$ and $\mathbf{u} \in \mathcal{M}_{R^m}$, $\mathbf{d} \in \mathcal{M}_{R^l}$, then the system is called forward complete. As usual, continuous function $\sigma: R_+ \to R_+$ belongs to class $\mathcal{K}$ if it is strictly increasing and $\sigma(0) = 0$; additionally it belongs to class $\mathcal{K}_\infty$ if it is also radially unbounded; and continuous function $\beta: R_+ \times R_+ \to R_+$ is from class $\mathcal{KL}$, if it is from class $\mathcal{K}$ for the first argument for any fixed second one, and it is strictly decreasing to zero by the second argument for any fixed first one.

Suppose that there exist functions $\rho_1, \chi_1 \in \mathcal{K}$ and constants $\rho_0, \chi_0 \in R_+$ such, that for all $\mathbf{x} \in R^n$ it holds:

$$|\mathbf{h}(\mathbf{x})| \leq \rho_1(|\eta(\mathbf{x})|) + \rho_0, \quad |\eta(\mathbf{x})| \leq \chi_1(|\mathbf{h}(\mathbf{x})|) + \chi_0. \tag{7}$$

**Definition 1**. *It is said that forward complete system (6) is practically state independent input-to-output stable (pSIIOS) with input norm operator $\mathcal{S}$ with respect to output $\mathbf{y}$ and input $\mathbf{d}$ (for $\mathbf{u}(t) \equiv 0$, $t \geq t_0 \geq 0$) if for all $\mathbf{x}(t_0) \in R^n$, $t_0 \geq 0$ and $\mathbf{d} \in \mathcal{M}_{R^l}$ there exist functions $\beta \in \mathcal{KL}$, $\gamma \in \mathcal{K}$ and $\sigma \in R_+$ such that*

$$|\mathbf{y}(t, \mathbf{x}(t_0), \mathbf{d})| \leq \beta(|\mathbf{h}(\mathbf{x}(t_0))|, t - t_0) + \gamma(\mathcal{S}[\mathbf{d}, t_0, t]) + \sigma$$

*for all $t \geq t_0$. The property pSIIOS satisfied for $\sigma = 0$ is called state independent input-to-output stability (SIIOS).* □

Possible choices of input norm operator $\mathcal{S}$ are $\mathcal{S}[\mathbf{d}, t_0, t] = \|\mathbf{d}\|_{[t_0, t)}$ (in this case for $\sigma = 0$ the property from definition 1 is reduced to well known SIIOS property from [29]) or integral one:

$$S[\mathbf{d},t_0,t] = \int_{t_0}^{t} \omega(|\mathbf{d}(\tau)|)\,d\tau,\ \omega \in \mathcal{K}.$$

The other closely connected input-output stability properties and relation between them can be found in [29], Lyapunov characterizations of these properties were presented in [30], small-gain theorem in [15]. In the case when $\mathbf{y} = \mathbf{x}$ the property is transforming to well known (practical) input-to-state stability property [28]. The following property is a local variant of SIIOS property for the case of inputs absence.

D e f i n i t i o n  2. *It is said that forward complete system (6) is locally uniformly asymptotically stable (LUAS) with respect to output $\psi$ if there exists constant $\Delta > 0$ such, that for all $\mathbf{x}(t_0) \in R^n$ with $|\eta(\mathbf{x}(t_0))| \leq \Delta$ and $\mathbf{u}(t) \equiv 0$, $\mathbf{d}(t) \equiv 0$, $t \geq t_0 \geq 0$ there exists function $\beta' \in \mathcal{KL}$ such that*

$$|\psi(t,\mathbf{x}(t_0),0)| \leq \beta'(|\eta(\mathbf{x}(t_0))|, t - t_0),\ t \geq t_0.\qquad\square$$

It is necessary to design control $\mathbf{u}: R^n \to R^m$ such, that
1) the outputs $\mathbf{y}$, $\psi$ are bounded in closed loop system for all initial conditions and inputs $\mathbf{d}$;
2) for the case $\mathbf{d}(t) \equiv 0$, $t \geq 0$ the following estimate holds for all $\mathbf{x}_0 \in R^n$:

$$|\psi(t,\mathbf{x}_0,0)| \leq \beta''(|\eta(\mathbf{x}_0)|, t),\ \beta'' \in \mathcal{KL}.$$

The last requirement is equivalent that for any $\lambda > 0$, $\kappa \geq 0$ and all $\mathbf{x}_0 \in R^n$, $|\eta(\mathbf{x}_0)| \leq \kappa$ there exists $T(\kappa,\lambda)$ such that

$$|\psi(t,\mathbf{x}_0,0)| \leq \lambda,\ \forall t \geq T(\kappa,\lambda),$$

function $T(\kappa,\lambda)$ is increasing with respect to the first argument for any fixed second one and it is strictly decreasing with respect to the second argument for any fixed first one. That implies finite time practical stabilization of the variable $\psi$ placed in the title of the paper. Let us introduce into consideration the following suppositions.

A s s u m p t i o n  1. *There exists "global" continuous control $\mathbf{u}_g : R^n \to R^m$ providing forward completeness property and pSIIOS property for the system (6) with input norm operator $S$ with respect to the output $\mathbf{y}$ and the input $\mathbf{d}$ for functions $\beta \in \mathcal{KL}$, $\gamma \in \mathcal{K}$.*$\qquad\square$

To design the control introduced in the assumption it is suggested to use CLF approach from the work [8], passification procedures [3], [26], [27] or integrator backstepping method [10].

A s s u m p t i o n  2. *There exists "local" continuous control $\mathbf{u}_l : R^n \to R^m$ providing LUAS property for output $\psi$ for $\beta' \in \mathcal{KL}$, $\Delta > 0$ and forward completeness property for the system (6), additionally for all $\mathbf{x}(t_0) \in R^n$, $t_0 \geq 0$ and $\mathbf{d} \in \mathcal{M}_{R^l}$ there exist functions $\alpha_1, \alpha_2, \alpha_3 \in \mathcal{K}$ and constant $\alpha \in R_+$ such that*:

$$|\psi(t,\mathbf{x}(t_0),\mathbf{d})| \leq \alpha_1(t - t_0) + \alpha_2(|\eta(\mathbf{x}(t_0))|) + \alpha_3(\|\mathbf{d}\|_{[t_0,t)}) + \alpha.\qquad\square$$

If the last estimate holds for $\psi = \mathbf{x}$, then it becomes an equivalent characterization of forward completeness property [2]. Therefore, the estimate additionally implies that for the control $\mathbf{u}_l$ the subsystem describing dynamics of the output $\psi$ is forward complete uniformly with respect to dynamics of all other components of the system state vector. The other necessary and sufficient conditions of forward completeness property can be also found in [2]. The following assumption together with (7) establishes a relations between the stabilized outputs $\mathbf{y}$ and $\psi$.

A s s u m p t i o n   3. *For the system (6) with the control* $\mathbf{u}_g$ *here exist* $\delta(\Delta) \in R_+$ *and* $T_{\Delta,\delta} > 0$ *such that for all* $\mathbf{x}_0 \in R^n$

$$|\mathbf{y}(t,\mathbf{x}_0,0)| \leq \delta(\Delta), \ \forall t \in [0, T_{\Delta,\delta}) \ \Rightarrow \ \exists t' \in [0, T_{\Delta,\delta}) : |\mathbf{\psi}(t',\mathbf{x}_0,0)| \leq \Delta.$$  □

The property from assumption 3 is a variant of detectability property (3) (zero-state-detectability from [4] or V-detectability from [24]), i.e.

$$\mathbf{y}(t) \equiv 0, \ t \geq 0 \ \Rightarrow \ \lim_{t \to +\infty} \mathbf{\psi}(t) = 0.$$

However there are two major differences here. At the first, the introduced property implies only that the set $|\mathbf{\psi}(t',\mathbf{x}_0,0)| \leq \Delta$ is reached at the time instant $t'$, but not convergence to this set (i.e. existence of a time instant $t'' > t'$ is possible such, that $|\mathbf{\psi}(t'',\mathbf{x}_0,0)| > \Delta$). At the second, the upper estimate $T_{\Delta,\delta}$ on the time for the desired set reaching is introduced in assumption 3. According to the relations (7) the property from assumption 3 trivially holds for $\Delta \geq \chi_0$.

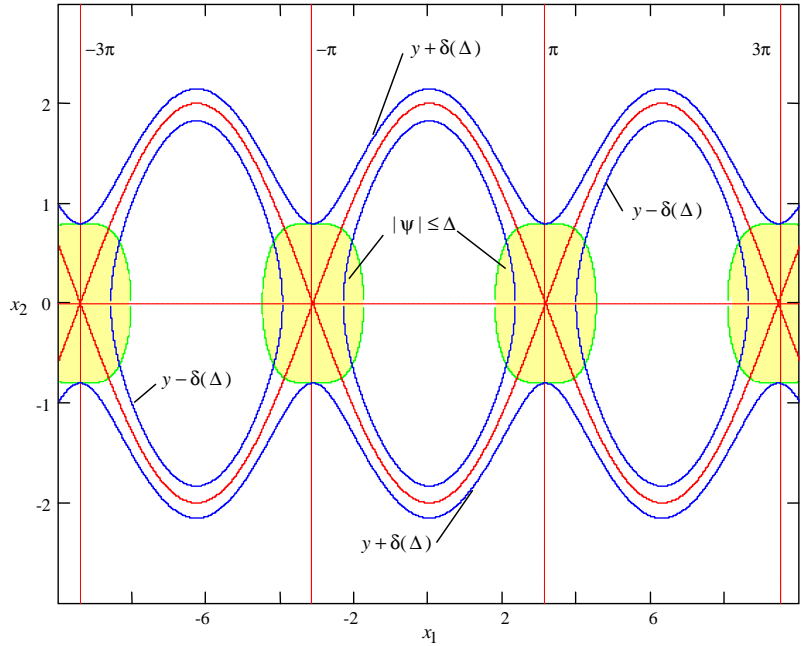

Fig. 3. Value of $\delta(\Delta)$ for pendulum (1).

For example, let us compute the values of $\delta(\Delta)$ and $T_{\Delta,\delta}$ for any given $\Delta \geq 0$ in the example considered in the previous section. In this case $y = H(\mathbf{x}) - H^* = 0.5 x_2^2 + \omega^2(1 - \cos(x_1)) - 2\omega^2$ and $\mathbf{\psi} = [1 + \cos(x_1), x_2]$ (these outputs have different dimensions). Then norms of the outputs have forms $|\mathbf{\psi}| = \sqrt{x_2^2 + (1+\cos(x_1))^2}$, $|y| = |0.5 x_2^2 - \omega^2(1+\cos(x_1))|$. Basing on these expressions (see also Fig. 3 for geometrical interpretation) it is possible to show that for any $\Delta \geq 0$ all trajectories staying in the set $|y| \leq \delta(\Delta)$ will intersect the set $|\mathbf{\psi}| \leq \Delta$ if $\delta(\Delta) = \min\{0.5\Delta^2, \omega^2 \Delta\}$. Indeed, all trajectories will intersect the set $|\mathbf{\psi}| \leq \Delta$ if boundaries of the set $|\mathbf{\psi}| \leq \Delta$ will intersect the boundaries of the set $|y| \leq \delta(\Delta)$ and the last set contains the first one. Since the set $|\mathbf{\psi}| \leq \Delta$ is an "ellipse" (the maximum deviations for the variables $x_1$, $x_2$ are achieved on the corresponding axes), the boundaries of the sets are intersecting on axes of the set $|\mathbf{\psi}| \leq \Delta$, where either

$x_1 = a\cos(\Delta - 1)$, $x_2 = 0$ or $x_1 = n\pi$, $x_2 = \Delta$. From this observation one can obtain desired expression for $\delta(\Delta)$. To calculate the time $T_{\Delta,\delta}$ simply note, that under such choice of $\delta(\Delta)$ trajectory should reach for the set $|\psi| \leq \Delta$ inside the set $|y| \leq \delta(\Delta)$ for the time less than a half of the system (1) time period $T_0$ (for the case without disturbances), thus in this case $T_{\Delta,\delta} = T_0/2$, and it is explicitly independent on values of $\delta(\Delta)$ and $\Delta$.

To solve the posed problem the following switching control is proposed

$$\mathbf{u}(t,\mathbf{x}) = \mathbf{u}_{i(t)}(\mathbf{x}), \tag{8}$$

where $i: R_+ \to \{l, g\}$ is the piecewise constant switching signal defining the type of control applied at the current time instant and assigned by the supervisor

$$t_{j+1} = \begin{cases} \arginf_{t \geq t_j + \tau_D} \mathbf{x}(t) \notin \mathbf{X}_\Delta & \text{if } i(t_j) = l; \\ \arginf_{t \geq t_j + \tau_D} \mathbf{x}(t) \in \mathbf{X}_l & \text{if } i(t_j) = g, \end{cases}$$

$$i(t_{j+1}) = \begin{cases} l, & \text{if } \mathbf{x}(t_{j+1}) \in \mathbf{X}_l; \\ g, & \text{if } \mathbf{x}(t_{j+1}) \notin \mathbf{X}_\Delta, \end{cases} \tag{9}$$

$$t_0 = 0, \ i(t_0) = \begin{cases} l, & \text{if } \mathbf{x}(t_0) \in \mathbf{X}_l; \\ g, & \text{otherwise,} \end{cases}$$

$$\mathbf{X}_l = \{\mathbf{x} : |\eta(\mathbf{x})| \leq \Delta\}, \ \mathbf{X}_\Delta = \{\mathbf{x} : |\eta(\mathbf{x})| \leq \beta'(\Delta, 0)\},$$

where $t_j$, $j = 1,2,3,...$ are the instants of switches, $j$ is the number of the last switch; $\tau_D > 0$ is the dwell-time constant. The control (8) equals to the control $\mathbf{u}_l$ into the set $\mathbf{X}_l$ and to the control $\mathbf{u}_g$ into the set $R^n \setminus \mathbf{X}_\Delta$. Signal $i(t)$ has constant value on the set $\mathcal{N} = \mathbf{X}_\Delta \setminus \mathbf{X}_l$. Since in generic case the set $\mathcal{N}$ can be non compact, it can not prevent fast switching phenomenon (chattering regime) rise. To bound the number of switches on any finite time interval the dwell-time constant $\tau_D$ is introduced. Such type of supervisors is called "dwell-time" one [16], [20], while the set of signal $i(t)$ constancy $\mathcal{N}$ plays a role of hysteresis in algorithm (9).

Theorem 1. *Let assumptions 1–3 hold. Then the control (8) with the supervisor (9) provide for the system (6) fulfillment of the following inequalities for all $\mathbf{x}_0 \in R^n$, $\mathbf{d} \in \mathcal{M}_{R^l}$ and $t \geq 0$:*

$$|\mathbf{y}(t,\mathbf{x}_0,\mathbf{d})| \leq \beta(\max\{|\mathbf{h}(\mathbf{x}_0)|, \rho_1(R_{\tau_D}^\Delta(\|\mathbf{d}\|)) + \rho_0\}, 0) + \gamma(S[\mathbf{d},0,t]) + \sigma, \tag{10}$$

$$|\psi(t,\mathbf{x}_0,\mathbf{d})| \leq \chi_1\left(3\beta(\rho_1(\max\{|\eta(\mathbf{x}_0)|, R_{\tau_D}^\Delta(\|\mathbf{d}\|)\}) + \rho_0, 0)\right) + \chi_1(3\gamma(S[\mathbf{d},0,t])) + \chi_1(3\sigma) + \chi_0, \tag{11}$$

*where $R_{\tau_D}^\Delta(\|\mathbf{d}\|) = \max\{\beta'(\Delta,0), \alpha_1(\tau_D) + \alpha_2(\Delta) + \alpha_3(\|\mathbf{d}\|) + \alpha\}$. If $\delta(\Delta) > \sigma$, then for the case $\mathbf{d}(t) \equiv 0$, $t \geq 0$ for all $\mathbf{x}_0 \in R^n$ for any $\lambda > 0$, $\kappa \geq 0$ with $|\eta(\mathbf{x}_0)| \leq \kappa$ there exists $T(\kappa,\lambda) \geq 0$ such that*

$$|\psi(t,\mathbf{x}_0,0)| \leq \lambda, \ t \geq T(\kappa,\lambda) = T_\kappa + T_{\Delta,\delta} + T_\lambda,$$

*where $\beta(\rho_1(\kappa) + \rho_0, T_\kappa) + \sigma \leq \delta(\Delta)$, $\beta'(\Delta, T_\lambda) \leq \lambda$.*

Proof. For the dwell-time supervisor algorithm on any time interval $[T_s, T_e)$ with $T_e > T_s \geq 0$ a finite number of

switches $N_{[T_s,T_e)}$ is possible and the following upper estimate holds:

$$N_{[T_s,T_e)} \leq 1+(T_e-T_s)\tau_D^{-1}.$$

Between switches combined control is continuous and equals to $\mathbf{u}_l$ or $\mathbf{u}_g$, which are continuous for all $\mathbf{x} \in R^n$. Thus resulting control is piecewise continuous function of time and solutions of the system (6), (8), (9) are continuous and defined for all $t \geq 0$. Indeed, on each interval of the control (8) continuity the system is forward complete (since the system (6) with the controls $\mathbf{u}_l$ or $\mathbf{u}_g$ are) and finite time escape phenomenon is not possible.

The interval of the system solutions existence $[0,+\infty)$ can be decomposer as follows:

$$[0,+\infty) = \Omega_l \cup \Omega_g,$$

where the control $\mathbf{u}_l$ is applied for all $t \in \Omega_l$ and the control $\mathbf{u}_g$ is used for $t \in \Omega_g$. From assumption 1 for all $t \in \Omega_g$, $\Omega_g = \bigcup_j [t_j,t_{j+1})$ (due to (9) each control has to be active during some finite time intervals with minimal length defined by $\tau_D$) the following estimate holds:

$$|\mathbf{y}(t,\mathbf{x}(t_j),\mathbf{d})| \leq \beta(|\mathbf{h}(\mathbf{x}(t_j))|,t-t_j) + \gamma(S[\mathbf{d},t_j,t]) + \sigma.$$

The possible variants for the values of $\mathbf{h}(\mathbf{x}(t_j))$ are as follows. At the first, $t_j = t_0$ and $\mathbf{h}(\mathbf{x}(t_j)) = \mathbf{h}(\mathbf{x}_0)$. Secondly, the control $\mathbf{u}_g$ can be switched on after the control $\mathbf{u}_l$ for $j > 1$. In this case two variants are possible. Firstly, the system leaves the set $\mathbf{X}_\Delta$ after dwell-time period, in this case $t_j > t_{j-1} + \tau_D$ and due to continuity of the solutions in this case

$$|\mathbf{\eta}(\mathbf{x}(t_j))| \leq \beta'(\Delta,0), \; |\mathbf{h}(\mathbf{x}(t_j))| \leq \rho_1(\beta'(\Delta,0)) + \rho_0.$$

In the second variant the system leaves the set $\mathbf{X}_\Delta$ for a time instant smaller than $\tau_D$ and in this case $t_j = t_{j-1} + \tau_D$. Recollecting forward completeness property estimate introduced in assumption 2 we obtain:

$$|\mathbf{\eta}(\mathbf{x}(t_j))| \leq \tilde{R}_{\tau_D}^\Delta(\|\mathbf{d}\|), \; |\mathbf{h}(\mathbf{x}(t_j))| \leq \rho_1(\tilde{R}_{\tau_D}^\Delta(\|\mathbf{d}\|)) + \rho_0,$$

where $\tilde{R}_{\tau_D}^\Delta(\|\mathbf{d}\|) = \alpha_1(\tau_D) + \alpha_2(\Delta) + \alpha_3(\|\mathbf{d}\|) + \alpha$. Thus, combining the estimates derived for both possible variants one can obtain for $j > 1$:

$$|\mathbf{\eta}(\mathbf{x}(t_j))| \leq R_{\tau_D}^\Delta(\|\mathbf{d}\|), \; |\mathbf{h}(\mathbf{x}(t_j))| \leq \rho_1(R_{\tau_D}^\Delta(\|\mathbf{d}\|)) + \rho_0.$$

By the same arguments the last estimates also hold for all $t \in \Omega_l$. Therefore, the estimate (10) is satisfied for all $t \geq 0$. Taking in mind (7) the estimate (10) implies (11).

Assume now that $\mathbf{d}(t) \equiv 0$, $t \geq 0$. Then according to assumption 1

$$|\mathbf{y}(t,\mathbf{x}_0,0)| \leq \beta(|\mathbf{h}(\mathbf{x}_0)|,t) + \sigma, \; t \geq 0$$

and since $\delta(\Delta) > \sigma$ there exists $T_\Delta \geq 0$ such, that

$$|\mathbf{y}(T_\Delta,\mathbf{x}_0,0)| \leq \beta(|\mathbf{h}(\mathbf{x}_0)|,T_\Delta) + \sigma = \delta(\Delta)$$

and for all $t \geq T_\Delta$ it holds that $|\mathbf{y}(t,\mathbf{x}_0,0)| \leq \delta(\Delta)$. Due to assumption 3 in this case

$$|\mathbf{\psi}(t',\mathbf{x}_0,0)| \leq \Delta$$

for some $t' \in [T_\Delta, T_\Delta + T_{\Delta,\delta})$. Further according to (9) the control law $\mathbf{u}_l$ will be switched on and from assumption 2 the following estimate holds:

$$|\psi(t,\mathbf{x}_0,0)| \leq \beta'(|\eta(\mathbf{x}(t'))|, t-t') \leq \beta'(\Delta, t-t'), \ t \geq t'$$

and the variable $\psi$ asymptotically converges to zero. For the variable $\psi$ for $t \in [0,t')$ according to the discussion above and (7) the following inequality holds for $t \in [0,t')$:

$$|\psi(t,\mathbf{x}_0,0)| \leq \chi_1(\beta(\rho_1(|\eta(\mathbf{x}_0)|) + \rho_0, t) + \sigma) + \chi_0.$$

Therefore, for arbitrary $\kappa \geq 0$ with $|\eta(\mathbf{x}_0)| \leq \kappa$ there exists $T_\kappa \geq 0$ such that

$$|\mathbf{y}(T_\kappa, \mathbf{x}_0, 0)| \leq \beta(\rho_1(\kappa) + \rho_0, T_\kappa) + \sigma \leq \delta(\Delta), \ t \geq T_\kappa$$

further for $t \geq T_\kappa + T_{\Delta,\delta}$ the LUAS estimate is satisfied from assumption 2 and for any $\lambda > 0$ there exists $T_\lambda \geq 0$ such that

$$|\psi(t,\mathbf{x}_0,0)| \leq \beta'(\Delta, T_\lambda) \leq \lambda, \ t \geq T(\kappa,\lambda) = T_\kappa + T_{\Delta,\delta} + T_\lambda.$$

Since both $T_\kappa$ and $T_\lambda$ are obtained for $\mathcal{KL}$ functions $\beta$ and $\beta'$ the function $T(\kappa,\lambda)$ is increasing in the first argument and decreasing in the second one. ∎

This theorem establishes conditions of the outputs $\mathbf{y}$ and $\psi$ boundedness for any $\mathbf{d} \in \mathcal{M}_{R^l}$ and convergence to zero of the variable $\psi$ in the noise-free case. These achieved results correspond to the control goals posed at the beginning.

The reason for the dwell-time constant $\tau_D$ introduction consists in the possible non compactness of set $\mathcal{N}$. If this is not the case, then in the supervisor (9) one can choose $\tau_D = 0$ (in this case the supervisor transforms to so-called "hysteresis" one [16], [20]) and theorem 1 admits the following modification.

T h e o r e m  2. *Let assumptions 1–3 hold and the set $\mathcal{N}$ be compact, then in (9) one can choose $\tau_D = 0$ and for the system (6), (8), (9) for all $\mathbf{x}_0 \in R^n$, $\mathbf{d} \in \mathcal{M}_{R^l}$ and $t \geq 0$ the following inequalities are satisfied:*

$$|\mathbf{y}(t,\mathbf{x}_0,\mathbf{d})| \leq \beta(\max\{|\mathbf{h}(\mathbf{x}_0)|, \rho_1(\beta'(\Delta, 0)) + \rho_0\}, 0) + \gamma(S[\mathbf{d}, 0, t]) + \sigma, \tag{12}$$

$$|\psi(t,\mathbf{x}_0,\mathbf{d})| \leq \chi_1\big(3\beta(\rho_1(\max\{|\eta(\mathbf{x}_0)|, \beta'(\Delta, 0)\}) + \rho_0, 0)\big) + \chi_1\big(3\gamma(S[\mathbf{d}, 0, t])\big) + \chi_1(3\sigma) + \chi_0. \tag{13}$$

*If $\delta(\Delta) > \sigma$, then for the case $\mathbf{d}(t) \equiv 0$, $t \geq 0$ for all $\mathbf{x}_0 \in R^n$ for any $\lambda > 0$, $\kappa \geq 0$ with $|\eta(\mathbf{x}_0)| \leq \kappa$ there exists $T(\kappa,\lambda) \geq 0$ such that $|\psi(t,\mathbf{x}_0,0)| \leq \lambda$ for all $t \geq T(\kappa,\lambda) = T_\kappa + T_{\Delta,\delta} + T_\lambda$, $\beta(\rho_1(\kappa) + \rho_0, T_\kappa) + \sigma \leq \delta(\Delta)$, $\beta'(\Delta, T_\lambda) \leq \lambda$.*

P r o o f. Let us introduce an upper bound for the right hand sides of the system (6) with the control (8) into the set $\mathcal{N}$:

$$F(r) = \max \left\{ \begin{array}{l} \sup_{\mathbf{x} \in \mathcal{N}, |\mathbf{d}| \leq r} \{|\mathbf{f}(\mathbf{x}, \mathbf{u}_l(\mathbf{x}), \mathbf{d})|\}, \\ \sup_{\mathbf{x} \in \mathcal{N}, |\mathbf{d}| \leq r} \{|\mathbf{f}(\mathbf{x}, \mathbf{u}_g(\mathbf{x}), \mathbf{d})|\} \end{array} \right\}.$$

The number $F(r)$ always exists and it is finite due to continuity property of the functions $\mathbf{f}$, $\mathbf{u}_l$, $\mathbf{u}_g$ and compactness of $\mathcal{N}$. Hence, the smallest time $\tau_r > 0$ that (6), (8) needs to cross set $\mathcal{N}$ can be estimated as $\tau_r = [\beta'(\Delta, 0) - \Delta] F(r)^{-1}$. On any time interval $[T_s, T_e)$ with $T_e > T_s \geq 0$ a finite number of switches $N^r_{[T_s, T_e)}$ is possible and the following upper estimate holds:

$$N^r_{[T_s,T_e)} \leq 1 + (T_e - T_s)\tau_r^{-1}.$$

Between switches combined control is continuous and equals to $\mathbf{u}_l$ or $\mathbf{u}_g$ controls, which are continuous for all $\mathbf{x} \in R^n$. Thus resulting control is piecewise continuous function of time and solutions of the system (6), (8), (9) are continuous and defined for all $t \geq 0$ (on each interval of the control (8) continuity the system is forward complete (since the system (6) with the controls $\mathbf{u}_l$ or $\mathbf{u}_g$ are) and finite time escape phenomenon is not possible).

The interval of the system solutions existence $[0,+\infty)$ can be decomposer as follows:

$$[0,+\infty) = \Omega_l \cup \Omega_g,$$

where the control $\mathbf{u}_l$ is applied for all $t \in \Omega_l$ and the control $\mathbf{u}_g$ is active for $t \in \Omega_g$. From assumption 1 for all $t \in \Omega_g$, $\Omega_g = \bigcup_j [t_j, t_{j+1})$ the following estimate holds:

$$|\mathbf{y}(t,\mathbf{x}(t_j),\mathbf{d})| \leq \beta(|\mathbf{h}(\mathbf{x}(t_j))|, t - t_j) + \gamma(S[\mathbf{d},t_j,t]) + \sigma.$$

The possible variants for values of $\mathbf{h}(\mathbf{x}(t_j))$ are as follows. At the first, $t_j = t_0$ and $\mathbf{h}(\mathbf{x}(t_j)) = \mathbf{h}(\mathbf{x}_0)$. At the second, the control $\mathbf{u}_g$ can be switched on after the control $\mathbf{u}_l$ and due to continuity of the solutions and (7) in this case

$$|\mathbf{\eta}(\mathbf{x}(t_j))| \leq \beta'(\Delta, 0), \quad |\mathbf{h}(\mathbf{x}(t_j))| \leq \rho_1(\beta'(\Delta,0)) + \rho_0.$$

For all $t \in \Omega_l$ the last estimates also hold. Therefore, the estimate (12) is satisfied for all $t \geq 0$. Tacking in mind (7) the estimate (13) follows from (12). Finite time convergence of the system to predefined neighbourhood of the set of zeros for output $\mathbf{\psi}$ can be proven applying the same arguments as in theorem 1. ∎

To apply the proposed in theorems 1 and 2 results in the motivating example, it is necessary to replace discontinuous "global" control (2) with its continuous approximation as the following one:

$$u = -u_m \varphi\left([H(\mathbf{x}) - H^*]x_2 \cos(x_1)\right), \tag{14}$$

$$\varphi(y) = \tanh(cy), \quad c > 0.$$

It is well known fact that the system (1) with the control (14) has additional equilibrium at the origin since in (14) $u(0) = 0$. In work [31] a discontinuous control was proposed, that ensures global stabilization of the upper equilibrium in finite time for almost all initial conditions (that is more it was proven, that it is not possible to provide global stabilization of the equilibrium via any continuous control). Under assumption that the system (6) with the control $\mathbf{u}_g$ is globally asymptotically stable for almost all initial conditions the theorems have the following corollary.

C o r o l l a r y   1. *Let assumptions 2,3 hold and there exist continuous control $\mathbf{u}_g : R^n \to R^m$ providing forward completeness property and for all $\mathbf{x}(t_0) \in R^n / \Xi$ (where $\Xi$ is a set of zero measure with property that if $\mathbf{x}(t_0) \in \Xi$, then $\mathbf{h}(\mathbf{x}(t,\mathbf{x}(t_0),0)) = \mathbf{h}(\mathbf{x}(t_0))$, $t \geq t_0$) for $\mathbf{d}(t) \equiv 0$, $t \geq t_0 \geq 0$ there exists function $\beta \in \mathcal{KL}$ such, that*

$$|\mathbf{y}(t,\mathbf{x}(t_0),0)| \leq \beta(|\mathbf{h}(\mathbf{x}(t_0))|, t - t_0), \quad t \geq t_0.$$

*Then for the case $\mathbf{d}(t) \equiv 0$, $t \geq 0$ the control (8) with the supervisor (9) provide for the system (6) fulfillment of the inequalities (10), (11) for all $\mathbf{x}_0 \in R^n$, $t \geq 0$. If $\delta(\Delta) > \sigma$, then for all $\mathbf{x}_0 \in R^n / \Xi$ for any $\lambda > 0$, $\kappa \geq 0$ with $|\mathbf{\eta}(\mathbf{x}_0)| \leq \kappa$*

*there exists* $T(\kappa, \lambda) \geq 0$:

$$|\psi(t, \mathbf{x}_0, 0)| \leq \lambda, \ t \geq T(\kappa, \lambda) = T_\kappa + T_{\Delta, \delta} + T_\lambda, \ \beta(\rho_1(\kappa) + \rho_0, T_\kappa) + \sigma \leq \delta(\Delta), \ \beta'(\Delta, T_\lambda) \leq \lambda.$$

P r o o f. The proofs of forward completeness and inequalities (10), (11) go similarly to the theorem 1 (taking into account that for $\mathbf{x}(t_0) \in \Xi$ it holds that $\mathbf{h}(\mathbf{x}(t, \mathbf{x}(t_0), 0)) = \mathbf{h}(\mathbf{x}(t_0))$, $t \geq t_0$). The proof of finite time convergence to predefined neighborhood of the output $\psi$ zeros also can be borrowed without modifications from theorem 1 for $\mathbf{x}_0 \in R^n / \Xi$. ∎

## IV. APPLICATION TO THE PENDULUM SYSTEM

Reconsider again the system (1). Let

$$\psi = [1 + \cos(x_1), x_2], \ y = 0.5 x_2^2 + \omega^2 (1 - \cos(x_1)) - 2\omega^2.$$

It was shown before how to calculate the value $\delta(\Delta)$ for arbitrary $\Delta > 0$ in this example with $T_{\Delta, \delta} = T_0 / 2$, thus assumption 3 holds. In the neighborhood of the points $(n\pi, 0)$, $n = \pm 1, \pm 3, \ldots$ the system (1) admits local approximation:

$$\begin{aligned} \dot{x}_1 &= x_2; \\ \dot{x}_2 &= -u + \vartheta, \end{aligned} \tag{15}$$

where $\vartheta \in R$ is an auxiliary disturbance. System (15) with the control

$$u_l = \alpha(k+1)[x_1 - n\pi + x_2], \ k > 0.5\omega^4 \tag{16}$$

for any positive $\alpha \geq 1$ admits locally the following Lyapunov function:

$$V(\mathbf{x}) = 0.5\omega^4 \left[ (x_1 - n\pi)^2 + (x_1 - n\pi + x_2)^2 \right], \ 0.25\omega^4 |(x_1 - n\pi, x_2)|^2 \leq V(\mathbf{x}) \leq 1.5\omega^4 |(x_1 - n\pi, x_2)|^2,$$

$$\dot{V} \leq -\omega^4 (x_1 - n\pi)^2 - \alpha\omega^4 \left( k + 1 - \alpha^{-1}(1 + 0.5\omega^4) \right)(x_1 - n\pi + x_2)^2 + 0.5\vartheta^2.$$

The system (1) in the set where $|\psi| \leq \Delta$, $\Delta < 1$ can be rewritten as a variant of (15) as follows:

$$\begin{aligned} \dot{x}_1 &= x_2; \\ \dot{x}_2 &= -bu + \vartheta, \end{aligned} \tag{17}$$

where $1 - \Delta \leq b \leq 1$ and $\vartheta = -\omega^2 \sin(x_1)$, $|\vartheta| \leq \omega^2 |x_1 - n\pi|$. For the system (17) and the control (16) with $\alpha = (1-\Delta)^{-1}$ the following estimate holds:

$$\dot{V} \leq -\kappa V, \ \kappa = \min\left\{ 1, 2(1-\Delta)^{-1} \left[ k + 1 - \frac{1 + 0.5\omega^4}{1 - \Delta} \right] \right\}.$$

Since locally $|1 + \cos(x_1)| \leq |x_1 - n\pi|$ and $|x_1 - n\pi| \leq 2\sqrt{1 + \cos(x_1)}$, then

$$|\psi(t)| \leq \sigma(|\psi(0)|)e^{-0.5\kappa t}, \ \sigma(s) = 2\sqrt{6} \max\{s, \sqrt{s}\}, \ \beta'(s, r) = \sigma(s)e^{-0.5\kappa r}.$$

Due to linear nature of the control the system is forward complete for any disturbance $\vartheta \in \mathcal{M}_R$. As it was just established, the output $\psi$ is proportional to the state vector, thus the estimate on the output growth is satisfied and assumption 2 holds.

The calculated function $\beta'$ defines the supervisor algorithm (9). The control (14) provides global stabilization of any specified level of the pendulum (1) energy for all initial condition except the origin [24] and all conditions of corollary 1 are satisfied. Then the global control is defined by (14), while the local control has the form (16) with $\alpha = (1-\Delta)^{-1}$. Let $\tau_D = 0$ in the supervisor (9). The results of the system simulation for the same initial conditions as in section 2 and for $\omega=1$, $u_m=0.1$, $c=20$, $k=1$ and $\Delta=0.2$ are shown in Fig. 4 and Fig. 5 and

$$E_x = 6.9, \quad E_H = 6.6.$$

The value of $E_x$ shows, that the system has non stochastic behavior with respect to the upper equilibrium convergence time. According to the result of corollary 1 the time shift between the convergence of energy (output $y$) and the stabilization of the upper equilibrium (output $\psi$) should be proportional to $T_{\Delta,\delta}$, that is confirmed by the results of simulation (for the chosen values of the system (1) parameters $T_0 = 7.41$), $T_\lambda = -2\kappa^{-1} \ln[\lambda / \sigma(\Delta)]$.

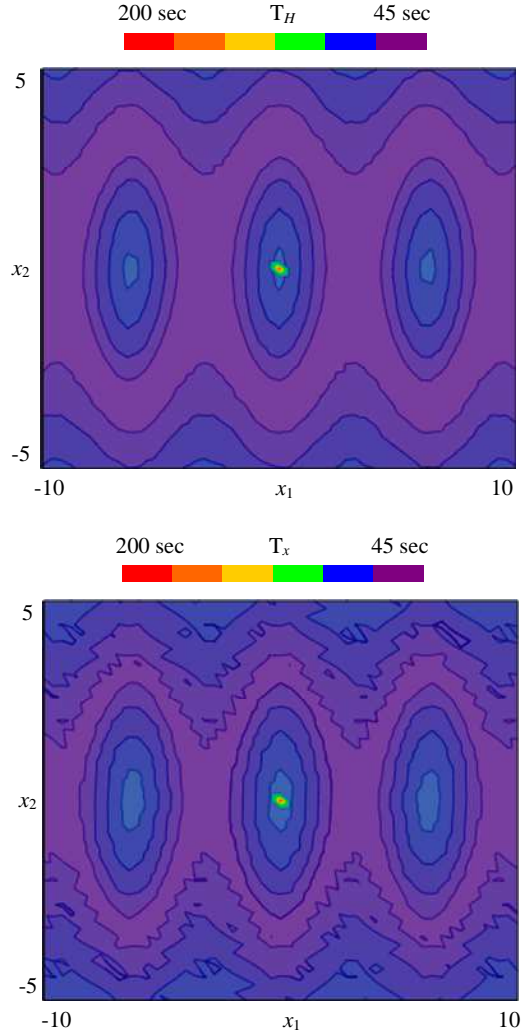

Fig. 4. Times $T_x$ and $T_H$ for system (1), (9), (14), (16).

## V. Conclusion

A new solution to the problem of nonlinear dynamical systems stabilization under zero-state-detectability assumption or its analogues is presented. The proposed solution ensures finite time practical stabilization of the system and it is based on uniting local and global controls. The global control provides boundedness of the system solutions and output convergence to zero, while local one ensures finite time convergence to a predefined goal set inside the zero dynamics set. Computer simulation demonstrates potentiality of the proposed solution for a pendulum system. For the pendulum it is shown that the problem of the finite time swinging up can be formulated as the problem of suitably defined entropy minimization. Useful upper estimates on nonlinear systems solutions under persistently excited gains are established in the appendix.

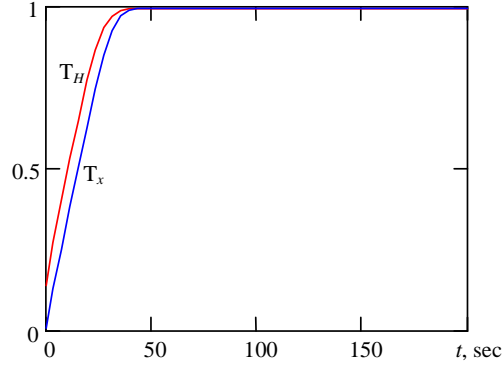

Fig. 5. Graphics of $T_x$ and $T_H$ for system (1), (9), (14), (16).

APPENDIX

D e f i n i t i o n   A 1 . *The Lebesgue measurable and square integrable matrix function* $\mathbf{R}: R_+ \to R^{l_1 \times l_2}$ *with dimension* $l_1 \times l_2$ *admits* $(L, \vartheta)$ *–persistency of excitation (PE) condition, if there exist strictly positive constants* $L$ *and* $\vartheta$ *such that for any* $t \geq 0$

$$\int_t^{t+L} \mathbf{R}(s)\mathbf{R}(s)^T \, ds \geq \vartheta \mathbf{I}_{l_1},$$

*where* $\mathbf{I}_{l_1}$ *denotes identity matrix of dimension* $l_1 \times l_1$.  □

The following lemma introduces an equivalent characterization of PE property used in the sequel.

L e m m a   A 1   [7]. *Let Lebesgue measurable and square integrable matrix function* $\mathbf{R}: R_+ \to R^{l_1 \times l_2}$ *with dimension* $l_1 \times l_2$ *be* $(L, \vartheta)$ *–PE. Then for any* $\ell \geq L$ *and* $t \geq 0$ *inequality is satisfied:*

$$\int_t^{t+\ell} \mathbf{R}(s)\mathbf{R}(s)^T \, ds \geq \frac{\vartheta}{2L} \ell \mathbf{I}_{l_1}.$$  ∎

The converse statement is obvious, if for matrix function $\mathbf{R}$ the last inequality is satisfied for all $\ell \geq L$, then $\mathbf{R}$ is $(L, \vartheta/2)$ –PE. Property from lemma A1 means that positive semidefinite matrix $\mathbf{R}(t)\mathbf{R}(t)^T$ has positive definite average matrix value for large enough time interval (the length of the interval should be bigger than $L$).

D e f i n i t i o n   A 2 . *The Lebesgue measurable and integrable matrix function* $\mathbf{A}: R_+ \to R^{l_1 \times l_1}$ *with dimension* $l_1 \times l_1$ *admits* $(L, \nu)$ *–positivity in average (PA) condition, if there exist constants* $L \geq 0$ *and* $\nu > 0$ *such, that for any* $t \geq 0$

$$\int_t^{t+\ell} \mathbf{A}(s) \, ds \geq \nu \ell \mathbf{I}_{l_1}.$$  □

The importance of PE or PA properties are explained in the following lemma.

L e m m a   A 2 . *Let us consider time-varying linear dynamical system*

$$\dot{\mathbf{p}} = -\mathbf{A}(t)\mathbf{p} + \mathbf{b}(t), \quad t_0 \geq 0, \tag{A1}$$

where $\mathbf{p} \in R^{l_1}$ and functions $\mathbf{A}: R_+ \to R^{l_1 \times l_1}$, $\mathbf{b}: R_+ \to R^{l_1}$ are Lebesgue measurable, $\mathbf{b}$ is essentially bounded, function $\mathbf{A}$ is $(L, \nu)$–PA for some $L > 0$, $\nu > 0$; $\mathbf{A}(t) = \mathbf{A}^T(t)$ and $\mathbf{A}(t) \geq -A\mathbf{I}_{l_1}$ for all $t \geq t_0$, $A \geq 0$. Then for any initial condition $\mathbf{p}(t_0) \in R^{l_1}$ solution of the system (A1) is defined for all $t \geq t_0$ and it admits the estimate

$$|\mathbf{p}(t)| \leq |\mathbf{p}(t_0)| e^{-\nu(t-t_0) - (1+\nu^{-1}A)L} + \|\mathbf{b}\| \{L + \nu^{-1} e^{-\nu L} + A^{-1}[e^{AL} - 1]\}.$$

P r o o f. Let us consider the expression for the system (A1) solutions:

$$\mathbf{p}(t) = \mathbf{p}(t_0) e^{-\int_{t_0}^{t} \mathbf{A}(\tau) d\tau} + \int_{t_0}^{t} e^{-\int_{\tau}^{t} \mathbf{A}(s) ds} \mathbf{b}(\tau) d\tau.$$

Matrix function $\mathbf{A}(t) \geq -A\mathbf{I}_{l_1}$ for all $t \geq t_0$ and $\|\mathbf{b}\| < +\infty$. Firstly, consider the case $A > 0$, then the inequality

$$|\mathbf{p}(t)| \leq |\mathbf{p}(t_0)| e^{A(t-t_0)} + \int_{t_0}^{t} \|\mathbf{b}\| e^{A(t-\tau)} d\tau \leq |\mathbf{p}(t_0)| e^{A(t-t_0)} + \|\mathbf{b}\| A^{-1}[e^{A(t-t_0)} - 1]$$

holds. So, solutions of the system (A1) are defined for all $t \geq t_0$. For $t \leq t_0 + L$ the last estimate takes form:

$$|\mathbf{p}(t)| \leq |\mathbf{p}(t_0)| e^{AL} + \|\mathbf{b}\| A^{-1}[e^{AL} - 1]$$

Using $(L, \nu)$–PA property of function $\mathbf{A}(t)$ we have for $t \geq t_0 + L$ and $\tau \leq t - L$:

$$-\int_{t_0}^{t} \mathbf{A}(\tau) d\tau \leq -\nu(t - t_0) \mathbf{I}_{l_1}, \quad -\int_{\tau}^{t} \mathbf{A}(s) ds \leq -\nu(t - \tau) \mathbf{I}_{l_1}.$$

Taking in mind that $|e^{\mathbf{S}}| \leq e^{\mu(\mathbf{S})}$ for a matrix $\mathbf{S} \in R^{n \times n}$, where $\mu(\mathbf{S}) = 0.5 \lambda_{\max}(\mathbf{S}^T + \mathbf{S})$ is the logarithmic norm of the matrix $\mathbf{S}$ and $\lambda_{\max}(\cdot)$ is the maximal eigenvalue of the corresponding matrix, and that

$$\mu\left(-\int_{t_0}^{t} \mathbf{A}(\tau) d\tau\right) = 0.5 \lambda_{\max}\left(-\int_{t_0}^{t} \mathbf{A}^T(\tau) d\tau - \int_{t_0}^{t} \mathbf{A}(\tau) d\tau\right) = \lambda_{\max}\left(-\int_{t_0}^{t} \mathbf{A}(\tau) d\tau\right) \leq -\nu(t - t_0), \quad \mu\left(-\int_{\tau}^{t} \mathbf{A}(\tau) ds\right) \leq -\nu(t - \tau),$$

we obtain for $t \geq t_0 + L$

$$|\mathbf{p}(t)| \leq |\mathbf{p}(t_0)| e^{-\nu(t-t_0)} + \int_{t_0}^{t} |\mathbf{b}(\tau)| e^{\varphi(t,\tau)} d\tau, \quad \varphi(t, \tau) = \begin{cases} A(t - \tau), & \text{if } \tau > t - L; \\ -\nu(t - \tau), & \text{if } \tau \leq t - L. \end{cases}$$

Performing the following simple transformations

$$\int_{t_0}^{t} |\mathbf{b}(\tau)| e^{\varphi(t,\tau)} d\tau = \int_{t_0}^{t-L} |\mathbf{b}(\tau)| e^{\varphi(t,\tau)} d\tau + \int_{t-L}^{t} |\mathbf{b}(\tau)| e^{\varphi(t,\tau)} d\tau \leq \int_{t_0}^{t-L} \|\mathbf{b}\| e^{-\nu(t-\tau)} d\tau +$$

$$\int_{t-L}^{t} \|\mathbf{b}\| e^{A(t-\tau)} d\tau \leq \nu^{-1} \|\mathbf{b}\| e^{-\nu L} + A^{-1} \|\mathbf{b}\| [e^{AL} - 1]$$

we obtain upper estimate for the system (A1) solutions for $t \geq t_0 + L$:

$$|\mathbf{p}(t)| \leq |\mathbf{p}(t_0)| e^{-\nu(t-t_0)} + \|\mathbf{b}\| \{\nu^{-1} e^{-\nu L} + A^{-1}[e^{AL} - 1]\}.$$

upper estimate for the case $A > 0$ can be obtained combining the last one with bounds, substantiated for $t \leq t_0 + L$:

$$|\mathbf{p}(t)| \leq |\mathbf{p}(t_0)| e^{-\nu(t-t_0-(1+\nu^{-1}A)L)} + \|\mathbf{b}\| \{\nu^{-1} e^{-\nu L} + A^{-1}[e^{AL}-1]\}.$$

Further, upper estimate of the lemma can be easily obtained, combining it with the estimate proposed in [7] for the case $A = 0$:

$$|\mathbf{p}(t)| \leq |\mathbf{p}(t_0)| e^{-\nu(t-t_0-L)} + (L + \nu^{-1} e^{-\nu L}) \|\mathbf{b}\|.\qquad\blacksquare$$

Lemma A2 establishes ISS-like property for linear PA system (A1) (ISS and iISS like estimate can be also found in [7] for PE system (A1), a variant of lemma A2 for scalar PA function $\mathbf{A}$ was proven in [9]). Let us extend this result for completely nonlinear system.

L e m m a   A 3 . *Consider a nonlinear dynamical system*

$$\frac{d\mathbf{z}(\tau)}{d\tau} = \mathbf{f}(\mathbf{z}(\tau), \mathbf{d}(\tau)), \ \tau \in R_+, \ \mathbf{z} \in R^n, \tag{A2}$$

*where the function* $\mathbf{f}: R^{n+m} \to R^n$ *ensures existence of the system solutions, which admits for any Lebesgue measurable and essentially bounded input* $\mathbf{d}: R_+ \to R^m$ *and initial condition* $\mathbf{z}_0 \in R^n$ *the estimate*

$$|\mathbf{z}(\tau, \mathbf{z}_0, \mathbf{d})| \leq \beta(|\mathbf{z}_0|, \tau) + \gamma(\|\mathbf{d}\|_{[0,\tau)}), \ \beta \in \mathcal{KL}, \ \gamma \in \mathcal{K}.$$

*Let Lebesgue measurable function* $a: R_+ \to R$ *be* $(L, \nu)$–*PA for some* $L > 0$, $\nu > 0$. *Then the system*

$$\frac{d\mathbf{p}(t)}{dt} = a(t)\mathbf{f}(\mathbf{p}(t), \mathbf{d}(t)), \ t \in R_+, \tag{A3}$$

*possesses the following estimate for* $t \geq L$:

$$|\mathbf{p}(t, \mathbf{z}_0, \mathbf{d})| \leq \beta(|\mathbf{z}_0|, \nu t) + \gamma(\|\mathbf{d}\|_{[0,t)}).$$

P r o o f . The systems (A2) and (A3) are introduced in the different time scales $\tau$ and $t$ correspondingly. Define the following relation between $\tau$ and $t$:

$$\tau = \int_0^t a(s)\,ds. \tag{A4}$$

Since $a$ is $(L, \nu)$–PA, then $\tau \geq \nu t$, $t \geq L$ and such time scales replacement is admissible for $\tau \geq \tau^* \geq L\nu$ (for time instants $t < L$ integral in the right hand side of (A4) may take negative values generically). In this case for $t \geq L$, $\tau \geq \tau^*$:

$$\mathbf{p}(t, \mathbf{z}_0, \mathbf{d}) \triangleq \mathbf{z}(\int_0^t a(s)\,ds, \mathbf{z}_0, \mathbf{d}) = \mathbf{z}(\tau, \mathbf{z}_0, \mathbf{d}), \ \hat{\mathbf{d}}(t) \triangleq \mathbf{d}(\int_0^t a(s)\,ds) = \mathbf{d}(\tau), \ d\tau = a(t)\,dt,$$

$$\frac{d\mathbf{p}(t)}{a(t)\,dt} = \frac{d\mathbf{z}(\tau)}{d\tau} = \mathbf{f}(\mathbf{z}(\tau), \mathbf{d}(\tau)) = \mathbf{f}(\mathbf{p}(t), \hat{\mathbf{d}}(t)),$$

where the last equation corresponds to (A3) for time instants when $a(t) \neq 0$. Note that for time instants when $a(t) = 0$ the established relation between solutions of the systems (A2) and (A3) also holds since in this case the system (A4) solution does not change its value. Substituting the proposed relation in the estimate for solution of (A2) we obtain the desired result.   $\blacksquare$

Lemma A3 establishes the property of ISS stability preservation under PA gain multiplication for nonlinear systems. Extension of LemmaA3 result for the cases of IOS or iISS properties is straightforward.